\documentclass[11pt]{amsart}

\usepackage{latexsym}
\usepackage{amsthm,amssymb,amsmath,amsopn}
\usepackage[all]{xypic}
\usepackage{pstricks}

\newtheorem{theorem}{Theorem}[section]
\newtheorem{lemma}[theorem]{Lemma}
\newtheorem{corollary}[theorem]{Corollary}
\newtheorem{proposition}[theorem]{Proposition}
\newtheorem{defn}[theorem]{Definition}
\newtheorem{remark}[theorem]{Remark}
\newtheorem{definition}[theorem]{Definition}
\newtheorem{example}[theorem]{Example}

\DeclareMathOperator{\sgn}{sgn}

\DeclareMathOperator{\cyc}{cyc}
\renewcommand{\neg}{\mathrm{neg}} 
\DeclareMathOperator{\desc}{desc}

\renewcommand{\epsilon}{\varepsilon}
\renewcommand{\theta}{\vartheta}

\newcommand{\vthinspace}{\hspace{0.5pt}}

\def\g{\gamma}

\def\d{\delta}

\def\s{\sigma}

\def\R{\mathbf{R}}
\def\E{\mathbb{E}}
\def\N{\mathbf{N}}

\def\g{\gamma}

\def\AA{{\mathcal A}}

\newcommand\shuf{
\setlength{\unitlength}{.4pt}
\begin{picture}(40,20)
\put(10,2){\line(1,0){20}}
\put(10,2){\line(0,1){10}}
\put(20,2){\line(0,1){10}}
\put(30,2){\line(0,1){10}}
\end{picture}
}

\newcommand\sshuf{\!
\setlength{\unitlength}{.3pt}
\begin{picture}(40,20)
\put(10,2){\line(1,0){20}}
\put(10,2){\line(0,1){10}}
\put(20,2){\line(0,1){10}}
\put(30,2){\line(0,1){10}}
\end{picture}\!
}

\begin{document}

\title[The moments of  L\'evy area]{A combinatorial method for calculating
the moments of  L\'evy area}

\author{Daniel Levin}
\address{Mathematical Institute, University of Oxford, 24-29 St Giles', Oxford
OX1 3LB, United Kingdom}
\email{levin@maths.ox.ac.uk}
\thanks{The first author is supported by the EPSRC Fellowship ``Partial
differential equations --- A rough
path approach'' GR/S18526/01}
\author{Mark Wildon}
\address{Department of Mathematics, University of Wales, Swansea, Singleton 
Park, Swansea SA2 8PP, United Kingdom}
\email{m.j.wildon@swansea.ac.uk}
\thanks{The second author is supported by EPSRC Grant EP/D054664/1}
\date{\today \newline
     } 

\begin{abstract}
We present a new way to compute the moments of the L\'evy 
area of a two-dimensional Brownian motion. Our approach uses
iterated integrals 
and combinatorial arguments involving the shuffle product. 
\end{abstract}

\maketitle

\section {Introduction}

In this paper we present a new approach to the problem of finding 
the moments of the signed area swept out by a two-dimensional
Brownian motion. This is a classical problem of great importance, 
originally solved by L{\'e}vy (see~\cite{Levy}).

We begin by explaining how these moments may
be defined. 
Given a piecewise smooth path $\gamma_t : [0,T] \rightarrow \R^2$ 
we may complete it to a loop  $\bar{\gamma}$
by closing it with the chord from $\gamma_T$ to $\gamma_0$. 
We may then define its \emph{signed area} to be
\[ \iint_{\R^2} n(\bar{\gamma},x) dx \]
where $n(\bar{\gamma},x)$ is the winding number of $\bar{\gamma}$
about the point $x \in \R^2$. 

\psset{xunit=0.8cm, yunit=0.8cm}
\begin{figure}[h]
\begin{pspicture}(8,6)    
	\psline{<-}(4,3)(8,6)
    \psline{-}(0,0)(4,3)
    {\pscurve{-}(0,0)(2.5,1)(4,5)(6,2)(8,6)} 
    \psline{->}(6,2)(6.05,2.0)
    \rput(-0.4,0){$\gamma_0$}
    \rput(8.35,6.0){$\gamma_T$}
    \rput(6.3,4.0){$+$}
    \rput(3.85,3.5){$-$}
    \rput(1.7,0.9){$+$}
    \rput(3.85,2.0){$0$}
    \rput(0.05,0){\circle*{0.1}}
    \rput(8.05,6){\circle*{0.1}}
\end{pspicture}
\caption{Contributions to the signed area of 
$\gamma : [0,T] \rightarrow \R^2$.}
\end{figure}
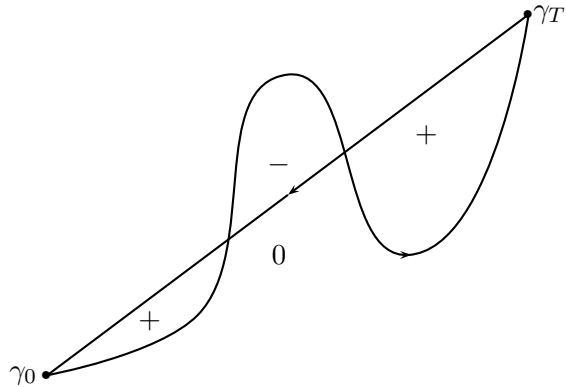

Now suppose that $B_t : [0,T] \rightarrow \R^2$ is
a two-dimensional Brownian motion. If  we complete~$B$ 
to a loop~$\bar{B}$ and attempt to define the 
signed area of~$B_t$ as before, then we immediately encounter the problem
that almost surely,~$n(\bar{B}_t, \cdot)$ is
not integrable on~$\R^2$ (see \cite[Theorem 55.I]{Levy3}). 
One solution is to replace~$B_t$ with
a sequence of piecewise linear dyadic approximations. 
In~\cite[Theorem 10]{Levy2} and~\cite[Chapter 55]{Levy3}, L{\'e}vy
proved that almost surely, the winding number integral
is defined for each approximation, and that the sequence
of areas converges. This gives one possible definition
of the L{\'e}vy area of the process~$B_t$. 

We may also define the signed
area corresponding to a smooth path  
$\gamma : [0,T] \rightarrow \R^2$ with $\gamma_t = (x_t,y_t)$
by
\[  \frac{1}{2} 
\int_0^T \bigl( (x_s-x_0) y'_s  - (y_s-x_0) x'_s \bigr) ds  
.\]
This observation motivates the following alternative definition of L{\'e}vy
area.

\begin{defn}\label{defn:levyarea}
Let $B_t = (X_t,Y_t)$ for $t\in [0,T]$ be a two-dimensional
Brownian motion starting at $0$. The \emph{L{\'e}vy area}
of $(B_t)_{0 \le t \le T}$ is given by the stochastic integral
\[ \AA_{T} = \frac{1}{2} \int_0^T (X_s dY_s - Y_s dX_s). \]
\end{defn}

\begin{remark}
In the sequel, we use $\AA$ as a shorter notation for $\AA_1$.
\end{remark}


L{\'e}vy showed in \cite{Levy} that almost surely the
definitions of signed area by dyadic approximation
and  by stochastic integration agree.
In his paper he also gave two different ways
to find the characteristic function, and hence
the moments, of $\AA_{T}$ when~$T = 2\pi$.

\begin{theorem}[L{\'e}vy]\label{thm:levy} If $T = 2\pi$ then
\begin{equation*}
\E \left( \exp \left(
i z \AA_{2\pi} \right) 
\right) =\left( \cosh \pi z 
\right)^{-1}.
\end{equation*}
\end{theorem}

L{\'e}vy's first proof uses the  definition of $\AA_{T}$ by
dyadic approximation.
His second starts from Definition~\ref{defn:levyarea}, 
but depends on earlier work by Kac, Siegert, Cameron and Martin (see
references cited in \cite[\S 1.6]{Levy}). 

In this paper we present
a direct and largely self-contained 
proof of Levy's Theorem, using Definition 1.1 to define L{\'e}vy area.
Our approach, which we outline in~\S 2 below,
is based on the fact that moments of L{\'e}vy area can be 
expressed as iterated integrals and hence calculated explicitly 
by exploiting the combinatorics of shuffle products. 

It seems likely that our methods can be applied
more broadly. We mention here that by 
using the multiplication
\[ (x,y,a)(x',y',a') = \Bigl( x+x',y+y',a+a' + \frac{1}{2}(xy'-yx') \Bigr) \]
we may identify points $(x,y,a) \in \R^3$ with elements
of the Heisenberg group. It is known (see \cite[Section 3.2.2]{LyLeCar}) 
that the process
$(X_t,Y_t,\AA_t)$ is a Brownian motion on this group.
Computing aspects of the joint distribution of $(X_t,Y_t,\AA_t)$
is a hard problem, involving Harish-Chandra formulae
(see \cite{Connell}); it is in effect aimed at understanding
the heat kernel on the Heisenberg group (see \cite{Gaveau}).
It seems likely that our approach may shed some light 
on these questions. Moreover, it should be possible to use
our methods to study the signed areas obtained when we replace
Brownian motion by measures related to higher order PDEs.

\section{Outline}
The outline of our proof is as follows. In \S 3 we 
use a simple scaling argument to show
that Theorem~\ref{thm:levy} is implied by the
following assertion about the moments of $\AA$.

\begin{theorem}\label{thm:moments}
$\E \AA^n = 2^{-n} E_n$.
\end{theorem}

\noindent Here $E_n$ is the $n$th Euler number, as 
defined by the generating
function
\[ \sum_{n=0}^\infty E_n \frac{z^n}{n!} = (\cos z)^{-1}.\]
The first few non-zero Euler numbers are $E_0 = 1$, $E_2 = 1$,
$E_4=5$, $E_6=61$. Of course all odd-numbered
Euler numbers are zero. (Correspondingly, one
can easily see that $\E \AA^n = 0$ if $n$ is odd.)
We then show that the moments of $\AA$ can be
expressed using iterated integrals.

In \S 4 we introduce the shuffle product on the
tensor algebra of a vector space, and use it to give 
an expression for $\E \AA^n$
as a certain coefficient in the
expansion of a shuffle product. In \S 5 we use a combinatorial
argument to determine this coefficient, thereby proving 
Theorem~\ref{thm:moments} and hence L{\'e}vy's
Theorem.




\section{The moments of L{\'e}vy's area for  Brownian motion}

We first show that L{\'e}vy's theorem (Theorem~\ref{thm:levy}) 
follows from
Theorem~\ref{thm:moments}. 
If we scale the Brownian path $B_t=(X_t, Y_t)$ defined for 
$0\le t \le 1$ by setting $\tilde{B}_s=\sqrt{T} B_{s/T}$
then we obtain a new Brownian path defined for
$0 \le t \le T$. As before, ~$\AA$ is the L\'evy area of~$B_t$ at time~$1$
and ~$\AA_{T}$ is the L\'evy area of $\tilde{B}_s$ at time $T$.
It follows easily from Definition~\ref{defn:levyarea} that~$\AA_{T}=T \AA_{1}$.
%
%
Hence, assuming that Theorem~\ref{thm:moments} holds, the moments of
L{\'e}vy area at time~$T$ are given by
\[ \E \AA_{T}^n = T^n n!\vthinspace E_n. \]
In particular, by setting $T = 2\pi$ we find that the
characteristic function of L{\'e}vy area at time $2\pi$ is 
\[
\E \left( \exp \left(i z \AA_{2\pi} \right) \right)=\sum_{n=0}^\infty \pi^n E_n
 \frac{(iz)^n}{n!} = (\cos \pi iz)^{-1}=(\cosh \pi z)^{-1},
\]
where we have absolute convergence of the series for $|z|<1/2$.

Therefore, to prove Theorem~\ref{thm:levy},
we may concentrate on finding the moments of $\AA$.
For this we shall need the following algebraic definition.


\begin{definition}
Let $V$ be a real vector space. Let
\[ 
T((V))=\prod_{k=0}^{\infty} V^{\otimes k}.
\]
Thus a typical element of $T((V))$ is a formal infinite sequence
$(a_0,a_1,a_2,\ldots)$ where $a_k \in V^{\otimes k}$. (By
convention $V^{\otimes 0} = \R$.) Clearly $T((V))$ is a real
vector space. It is easy to check that $T((V))$ becomes an
algebra with unit if we define the product of
$\mathbf{a}=(a_0, a_1, a_2, \dots),
\mathbf{b}=(b_0, b_1, b_2, \dots) \in T((V))$ by
\[ \mathbf{a} \otimes \mathbf{b} =\left(\dots, \sum_{j=0}^k a_j \otimes
b_{k-j}, \dots \right). \]
We also define the \emph{exponential} of a tensor $\mathbf{a}$
to be the formal sum
\[
\exp\left( {\bf a} \right)=\sum_{n=0}^{\infty} \frac{{\bf a}^{\otimes n}}{n!}.
\]
\end{definition}




We can now define  the signature of a Brownian
motion. As motivation, we first recall this definition for a path. 
If $\gamma : [0,T] \rightarrow \R^2$ 
is a path of finite length (see \cite{LyHa}) then its \emph{signature} 
is the formal infinite sum  
\[ \textbf{X}_{s,t}(\gamma) = \sum_{k=0}^{\infty} \int_{s<t_1<\dots<t_k<t} 
d\gamma_{t_1} \otimes \dots \otimes d\gamma_{t_k} \in T((\R^2)) \]
defined for $0 \le s < t \le T$.

\begin{example}
For example, if $e_0, e \in \R^2$ and $\gamma_t = e_0 + t e$, 
then
\[ \textbf{X}_{s,t}(\gamma) = \exp  (t-s) e . \]
\end{example}

One of the most important properties of the signature is that it is 
multiplicative (for the proof of the following theorem see \cite{Chen}). 

\begin{theorem}
For any $0\le r<s<t\le T$ 
\begin{equation*}
{\bf X}_{r,t}(\g)={\bf X}_{r,s}(\g) \otimes {\bf X}_{s,t}(\g).
\end{equation*}
\end{theorem}

\begin{definition}
Let $B : [0,T] \rightarrow \R^2$ be a two-dimensional 
Brownian motion. For $0 \le s <t \le T$ we define 
the \emph{signature} of $B$ to be the formal infinite sum
\begin{equation*}
    \textbf{X}_{s,t}(B) = \sum_{k=0}^{\infty} \int_{s<t_1<\ldots<t_k<t} 
    dB_{t_1}\otimes \dots \otimes dB_{t_k}.
\end{equation*}
\end{definition}

\begin{remark}
Up till now, the signature of a Brownian motion and other auxiliary objects 
could be interpreted either in the Stratonovich or It\^o sense. However, the 
next theorem holds only for the Stratonovich integral. For its
proof see either Fawcett \cite{Fawcett} or Lyons--Victoir \cite{LyVic}.
\end{remark}

\begin{theorem}
\label{LV}
If $B$ is a Brownian motion in $\R^2$ then
\begin{equation} 
\label{LV_result}
\E \bigl(\mathbf{X}_{0,1}  (B) \bigr)=\exp \left( \frac{1}{2} \left( 
e_1 \otimes e_1 + e_2 \otimes e_2 \right) \right)
\end{equation}
where $e_1$ and $e_2$ are any two orthogonal vectors in $\R^2$. 
\end{theorem}




In particular,~\eqref{LV_result} implies that 
\begin{equation}
\label{trunc_ID} 
\E \left(X^{2n}_{0,1} \left( B \right) \right)=\frac{1}{2^n n!} 
\left( e_1 \otimes e_1 + e_2 \otimes e_2 \right)^{\otimes n}
\end{equation}
where $X^{2n}_{0,1}$ is the component of $\mathbf{X}_{0,1}$
lying in $(\R^2)^{\otimes 2n}$.


\section {Shuffle products and other combinatorial objects}

We now introduce an important combinatorial object which will be used in the 
sequel.
\begin{definition}
\label{shuffle}
We define the set $S_{m,n}$ of $(m,n)$ shuffles to be the subset of 
permutations in the symmetric group $S_{m+n}$ defined by
\begin{equation*}
S_{m,n}=\{ \s \in S_{m+n}: \s(1) < \dots < \s(m), \; \s(m+1)<\dots<\s(m+n)\}.
\end{equation*}
\end{definition}

\begin{remark}
The term ``shuffle'' is used because such permutations arise when
one riffle 
shuffles a deck of $m+n$ cards cut into one pile of $m$ cards and a second
pile of $n$ cards.
\end{remark}

Let $V=\R^2$ with the orthogonal basis $e_1$ and $e_2$.  
Let $V^{*}$ be the space dual to $V$ and let $e^1$ and $e^2$ 
be its dual basis. 
Let $n \in \mathbf{N}$. The elements $e_{i_1} \otimes \dots \otimes 
e_{i_n}$, where each $i_k \in \{1,2\}$, for
$k=1, \dots, n$, form a basis of $V^{\otimes n}$. The
corresponding dual basis of ${\left( V^{*} \right)}^{\otimes n}$ is 
given by the elements $e^{i_1} \otimes 
\dots \otimes e^{i_n}$.
  
There is a natural duality $\left<\;, \; \right> : V^{\otimes n} \times
(V^{\star})^{\otimes n} \rightarrow \R$ defined by:
\begin{equation*}
\left<e^{i_1} \otimes \dots \otimes e^{i_n}, e_{j_1} \otimes \dots \otimes
e_{j_n}\right>=\d_{i_1 j_1} \dots \d_{i_n j_n}.
\end{equation*}

\begin{definition}
Set $(k_1, \dots k_{m+n})=(i_1, \dots, i_m, j_1, \dots, j_n)$. The shuffle 
product of two tensors
$e^I=e^{i_1}\otimes \dots \otimes e^{i_m}$ and 
$e^J=e^{j_1}\otimes \dots \otimes e^{j_n}$, is the tensor 
$e^I \shuf e^J$ defined by
\begin{equation*}
e^I \shuf e^J=\sum_{\s \in S_{m,n}} e^{k_{\s^{-1}(1)}} \otimes \dots \otimes 
e^{k_{\s^{-1}(m+n)}}. 
\end{equation*}
Let
$T(V^\star) = \bigoplus_{k=0}^\infty (V^\star)^{\otimes k}$ be the ordinary tensor algebra on $V^\star$.
The shuffle product~$\shuf$ extends to a bilinear map 
\[ T(V^\star) \times T(V^\star) \rightarrow T(V^\star).\]
\end{definition}

For example, the reader may check that
\[  (e \otimes f) \shuf g = e \otimes f \otimes g + e \otimes g \otimes f
+ g \otimes e \otimes f \]
for any $e,f,g \in (\R^2)^\star$. For an alternative description
of the shuffle product, see Definition~5.14 in \cite{LyHa}.

\begin{remark}
It follows easily from the definition that the
shuffle product is commutative and 
associative. We shall use the following notation for the shuffle product 
applied $N$ times:
\begin{equation*}
{\bf a}^{\sshuf N}=\underbrace{{\bf a} \shuf \dots \shuf {\bf a}}_N.
\end{equation*}
Later in \S 5 we shall also use the analogous version of the shuffle
product defined on the tensor powers of $V$.
\end{remark}

For each path $\g_s$, $s \in [0,T]$ of finite
length we now introduce a real-valued function
\[ \varphi(\g): T(V^{\star}) \rightarrow \R \] 
defined on the tensor 
${\bf e}=e^{i_1} \otimes \dots \otimes e^{i_n}$ by  
\begin{align}
\varphi_{\bf e} (\g)
\notag 
&= \left<e^{i_1} \otimes \dots \otimes e^{i_n}, 
\int_{0<t_1<\dots<t_n<T} d\g_{t_1} \otimes \dots \otimes d\g_{t_n}\right>
\notag \\
&=\int_{0<t_1<\dots<t_n<T} d\g^{i_1}_{t_1} \dots d\g^{i_n}_{t_n}.
\label{coord}
\end{align}
The fundamental property of 
$\varphi_{\bf e}(\g)$ (see \cite[Theorem 2.15]{LyLeCar}) is that   
\begin{equation}
\label{shuffle_prod}
\varphi_{\bf e}(\g) \varphi_{\bf f}(\g)=\varphi_{{\bf e} \sshuf {\bf f}}(\g).
\end{equation}
From now on we change slightly notations for $e_1, e_2$ and $e^1, e^2$. Let 
$x$, $y$ be a basis for $\R^2$ and let
 $x^{*}$, $y^{*}$ be the dual basis of $(\R^2)^\star$. 
We use these techniques together to prove the following theorem.
\begin{theorem}
\label{Main}
The $n$th moment of L\'evy area at time $1$ for a two-dimensional
Brownian motion~$B$
starting at zero is the signature 
${\bf X}_{0,1}(B)$ contracted with the shuffle powers of the dual 
tensor $\frac{1}{2}\left( x^{*} \otimes y^{*} - y^{*} \otimes x^{*} \right)$.
That is,
\begin{equation*}
\E \left( \AA^n \right)=
2^{-2n} \left<(x^{*} \otimes y^{*} - y^{*} \otimes x^{*})^{\sshuf n}, 
\frac{(x \otimes x + y \otimes y)^{\otimes n}}{n!}\right>.
\end{equation*}
\end{theorem}

\begin{proof}
Let $B=(X_s,Y_s)$, $0\le s \le 1$ be a Brownian path in $\R^2$ starting at 
zero. Then 
\begin{equation*}
\AA=\frac{1}{2} \int_0^1 \left( X_s \, d Y_s - Y_s \, d X_s \right)=
\frac{1}{2} \int_{0<t<s<1} \left( d X_{t} \, d Y_{s} - d Y_t \, d X_s \right).
\end{equation*} 
In a more canonical notation as in \eqref{coord} we may write this down as 
\begin{equation*}
\AA=\frac{1}{2} \int_{0<t_1<t_2<1} \left( d X_{t_1} \, d Y_{t_2} - d Y_{t_1} 
\, d X_{t_2} \right).
\end{equation*}
So $\AA=\varphi_{\frac{1}{2}\left( x^{*} \otimes y^{*} - y^{*} \otimes x^{*} 
\right)} 
\left( B \right)=\frac{1}{2} \left( \varphi_{x^{*} \otimes y^{*}} \left(B 
\right) - \varphi_{y^{*} \otimes x^{*}} \left( B \right) \right)$. 
Further, using~(\ref{shuffle_prod}) we have 
\begin{align*}
\AA^n&=\frac{1}{2^n} \Bigl<\left( x^{*}  \otimes y^{*} - y^{*} \otimes x^{*} 
\right), \int_{0<t_1<t_2<1} dB_{t_1} \otimes dB_{t_2}\Bigr>^n
\\
&=\frac{1}{2^n} \Bigl<\left( x^{*}  \otimes y^{*} - y^{*} \otimes x^{*} 
\right)^{\sshuf n}, \;
X_{0,1}^{2n}\left( B \right)\Bigr>.
\end{align*} 

Taking the expectation of $\AA^n$ we have by \eqref{trunc_ID},
\begin{align*}
\E\left( \AA^n \right) &= 2^{-n} \left<\left( x^{*} \otimes y^{*} - y^{*} 
\otimes x^{*} \right)^{\sshuf n}, \; \E \left( X_{0,1}^{2n}\left(B \right)
\right)\right>
\\ 
&=2^{-2n} \left<\left( x^{*} \otimes y^{*} - y^{*} \otimes x^{*} 
\right)^{\sshuf n}, \; 
\frac{\left( x \otimes x + y \otimes y \right)^{\otimes n}}{n!}\right>. 
\end{align*}

\vspace{-22pt}

\end{proof} 

Hence to prove Theorem~\ref{thm:moments} it is sufficient to prove
the following theorem.

\begin{theorem}\label{thm:MW1} 
\[ \left<\left( x^{*} \otimes y^{*} - y^{*} \otimes x^{*} \right)^{\sshuf n}, 
\; 
\left( x \otimes x + y \otimes y \right)^{\otimes n}\right>
= 2^n n!\vthinspace E_n.
\]
\end{theorem}

\section{Proof of Theorem~\ref{thm:MW1}}  
\label{Proof_Theor}

\subsection{}
From now on it will often be convenient to use
a shorter notation for elements in the standard basis
of the tensor
algebra~$T(V)$, in which we write 
$xy$ rather than~\hbox{$x\otimes y$},~$x^2$ rather than~\hbox{$x
\otimes x$},
and so on. Using this notation 
the standard basis elements of~$T(V)$ are simply the words
in the letters~$x$ and~$y$. 

We shall say that a 
word in the letters $x$ and $y$ is \emph{even} if (i) it is of the
form 
$z_1^2 \ldots z_n^2$ where each $z_i \in \{x,y\}$, and (ii)
there are equal numbers of $x$'s and $y$'s. 
When we expand
\[ (x \otimes x + y \otimes y)^{\otimes n} \]
we obtain the sum of all words in $x$ and $y$ of length $2n$
satisfying condition~(i). If such a word is not killed
by $(x^\star \otimes y^\star - y^\star \otimes x^\star)^{\sshuf n}$ then
clearly it must also satisfy~(ii). Hence Theorem~\ref{thm:MW1} is equivalent
to the following assertion.

\begin{theorem}\label{thm:eqvform}
Let 
\[ (xy-yx)^{\sshuf n} = \sum_v \lambda_v v \]
where the sum is over all words $v$ of length $2n$. Let
$u_n = \sum \lambda_w$ 
where the sum is over all even words $w$ of length $2n$. Then
\begin{equation*} u_n = 2^n n!\vthinspace E_n 
\end{equation*}
\end{theorem}


\subsection{}
We begin the proof of 
Theorem~\ref{thm:eqvform} by noting that if $n$ is odd then there are
no even words of length $2n$, and so $u_n = 0$, as required.
We may therefore assume that~$n = 2m$ is even, so the
even words that appear in $(xy-yx)^{\sshuf n}$ are of length $4m$,
with the pairs $xx$ and $yy$ each appearing exactly $m$ times. 

Our proof
depends on counting the combinatorial objects
introduced in the next definition.

\begin{defn}An \emph{$xy$-matching} is a pair $(w,\sigma)$
where $w$ is an even word, of length $2m$ say,
and $\sigma$ is a fixed-point-free involution in
the symmetric group $S_{2m}$
such that $w_i = x$ if and only if $w_{\sigma(i)} = y$.

Given an $xy$-matching $\delta = (w,\sigma)$
we define the \emph{negativity} of $\delta$ by
\[ \neg(\delta) = \# \{i \in \{1 \ldots 2m\} : 
\text{$w_i = x$ and $\sigma(i) < i$}\}.
\]
We define the \emph{sign} of $\delta$ by $\sgn(\delta) = (-1)^{\neg(\delta)}$.
Let $N_t(w)$ be the number of $xy$-matchings with underlying word
$w$ and negativity $t$.
\end{defn}

It will be very useful to represent $xy$-matchings by diagrams
such as the one below.

\bigskip
\bigskip
\hfil
\xymatrix@C=6pt{x \ar@{-}@/^12pt/[rr]&
x \ar@{-}@/^24pt/@<2pt>[rrrrrr]&
y &
y \ar@{-}@/_12pt/[rr]&
x \ar@{-}@/^12pt/[rr]&
x &
y &
y &
}
\hfill

\bigskip
\begin{center}
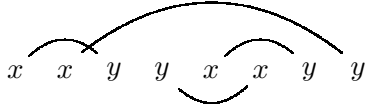
{\sc Figure 2.} The $xy$-matching $(xxyyxxyy, (13)(28)(46)(57))$.
\end{center}

\bigskip

\noindent The arcs contributing to $\neg(\delta)$ are those drawn
below the word, thus here $\neg (\delta) = 1$ and
$\sgn(\delta) = -1$.
As an exercise, the reader may check that
there are in total $16$ $xy$-matchings
with underlying word $xxyyxxyy$ and negativity~$1$, and so
$N_1(xxyyxxyy) = 16$.


\begin{proposition}\label{prop:nasty}
Let $w$ be an even word of length $2m$ and
let $s + t = 2m$. The coefficient
of $w$ in $(xy)^{\sshuf s} \shuf (yx)^{\sshuf t}$ is
$s!t!\vthinspace N_t(w)$.
\end{proposition}

The proof of this proposition
is postponed to the appendix at the end of this 
paper. The idea is to 
associate to each $xy$-matching
with underlying word~$w$ and negativity~$t$ exactly
$s!t!$ ways to obtain $w$ by expanding the shuffle product
\hbox{$(xy)^{\sshuf s} \shuf (yx)^{\sshuf t}$}. With the help
of a formal definition of an expansion of a shuffle product
of this form,
we are able to
show that these possibilities are exhaustive.
Here we shall illustrate
the correspondence when
$m = 4$, $w = xxyyxxyy$,~$s=3$ and~$t=1$.

The figure below shows one
way to obtain $w$ by expanding the shuffle prouduct~$xy \shuf xy
\shuf xy \shuf yx$.

\medskip
\newcommand{\y}{\raisebox{3pt}[10pt]{$y$}}
\vbox{
\hfil
\xymatrix@C=3pt@R=42pt{
x\ar@{.>}[d] & \y\ar@{.>}[drr] 
& \shuf & x\ar@{.>}[dll] & \y\ar@{.>}[drrrrrr]  & \shuf 
& x\ar@{.>}[d] & \y\ar@{.>}[drr]  &  \shuf & \y\ar@{-->}[dlllll] & x\ar@{-->}[dlll]  \\
x & x & 
& y & y& & x  & x & & y & y}
\hfill

\smallskip
\begin{center}
{\sc Figure 3.} An expansion of $xy \shuf xy \shuf xy \shuf yx$.
\end{center}}

\medskip
\noindent

We obtain the corresponding $xy$-matching by connecting the letters
coming from the same $xy$ or $yx$ term in $xy \shuf xy \shuf xy \shuf yx$.

\bigskip\bigskip
\hfil
\newcommand{\cupspace}{\hspace{8pt}}
\xymatrix@C=3pt@R=36pt{
x\ar@{-}@/^12pt/[rrr] & x \ar@{-}@/^24pt/[rrrrrrrrr] & \cupspace
& y & y\ar@{-}@/_12pt/[rrr] & \cupspace& x\ar@{-}@/^12pt/[rrr] & x 
& \cupspace& y & y}
\hfill

\bigskip
\begin{center}
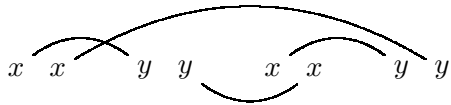
{\sc Figure 4.} The $xy$-matching corresponding to the expansion in
Figure 3.
\end{center}

\medskip
\noindent 
Note that this matching, $(xxyyxxyy, (13)(28)(46)(57))$,
has negativity~$1$, corresponding to the single $yx$ term.

There are in total $3! 1!$
ways to obtain~$w$ by expanding the shuffle product 
$xy \shuf xy
\shuf xy \shuf yx$
which correspond to this matching. The
remaining five are obtained by permuting identical words
in the top line of Figure~3.  For example,
we could just as well get the loop joining the first
$x$ to the first $y$ 
by expanding the shuffle product as shown below.

\medskip
\hfil
\xymatrix@C=3pt@R=42pt{
x\ar@{.>}[dr] & \y\ar@{.>}[drrrrrrrrr] 
& \shuf & x\ar@{.>}[dlll] & \y\ar@{.>}[dl]  & \shuf 
& x\ar@{.>}[d] & \y\ar@{.>}[drr]  
&  \shuf & \y\ar@{-->}[dlllll] & x\ar@{-->}[dlll]  \\
x\ar@{-}@/^12pt/@<-2pt>[rrr] & x \ar@{-}@/^18pt/[rrrrrrrrr] & \cupspace
& y & y\ar@{-}@/_12pt/[rrr] & \cupspace& x\ar@{-}@/^12pt/@<-1pt>[rrr] & x 
& \cupspace & y & y}
\hfill

\medskip
\begin{center}
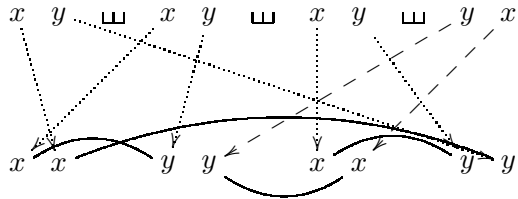
{\sc Figure 5.} Another of the $3!1!$ expansions of $xy \shuf xy \shuf xy \shuf yx$
giving the $xy$-matching shown in Figure~4.
\end{center}


\begin{corollary}\label{cor:nasty} Let $w$ be an even word of length $2m$.
The coefficient of $w$ in $(xy-yx)^{\sshuf 2m}$ is
\[ (2m)! \sum_{\delta} \sgn(\delta) \]
where the sum is over all $xy$-matchings with underlying word $w$.
\end{corollary}

\begin{proof}Let $s + t = 2m$. As the shuffle
product is commutative,
when we expand $(xy-yx)^{\sshuf 2m}$ we obtain
$(-1)^t (xy)^{\sshuf s} \shuf (yx)^{\sshuf t}$ exactly $\binom{2m}{s}$
times.
Hence,
by Proposition~\ref{prop:nasty}, the coefficient we seek is
\[ \sum_{s+t = 2m} \binom{2m}{s} (-1)^t s!t! N_t(w) =
(2m)! \sum_{t=0}^{2m} (-1)^t N_t(w) 
= (2m)! \sum_{\delta} \sgn(\delta). \]
where the final equality holds because
the middle sum counts each $xy$-matching exactly
once, with the appropriate sign.
\end{proof}

\subsection{}
We now introduce our second and final
combinatorial object.

\begin{defn}An \emph{$XY$-matching} is a pair $(W,\sigma)$
where $W$ is a word in letters $X$ and $Y$ with equal
numbers of $X\!$'s and $Y\!$'s, say $m$ of each, and 
$\sigma \in S_{m}$ is a permutation such that
$W_i = X$ if and only if $W_{\sigma(i)} = Y$. 
\end{defn}

Once again, it is very useful to represent
$XY$-matchings  by diagrams. The figure below
shows a typical example.

\vbox{
\bigskip
\bigskip
\bigskip

\hfil
\xymatrix{A \ar@/^12pt/[r]& 
B \ar@/_12pt/[r]& 
A \ar@/^12pt/[r]& 
B \ar@/_15pt/@<-9pt>[lll]& 
B \ar@/_9pt/[r]& 
A \ar@/^15pt/[l]
}
\hfill

\bigskip
\begin{center}
{\sc Figure 5.} The $XY$-matching $(XYXYYX, (1234)(56))$.
\end{center}}
\medskip

As before, we need certain quantities associated with
an $XY$-matching.

\begin{defn}
Let $\Delta = (W,\sigma)$ be an $XY$-matching. The \emph{length} of $\Delta$
is the length of the word $W$. 
We define the \emph{negativity} of $\Delta$ by
\[ \neg (\Delta) = \# \{ i : \text{$W_i = X$ and $\sigma(i) < i$ \emph{or}
$W_i = Y$ and $\sigma(i) > i$} \}. \]
We define the \emph{sign} of $\Delta$ by
\[ \sgn(\Delta) = (-1)^{\neg(\Delta)}.\] 
We define the \emph{cycle count} of $\Delta$ by
\[ \cyc(\Delta) = \text{$\#$disjoint cycles in the permutation $\sigma$}.\]
Finally we let $e(\Delta)$
be the even word obtained from $\Delta$ by replacing each $X$
with $xx$ and each $Y$ with $yy$.
\end{defn}

For example, if $\Delta$ is the $XY$-matching shown in Figure 5 above, 
then $\Delta$ has length~$6$,
$\neg(\Delta) = 3$ (the $3$ arcs drawn below the word contributing),
$\sgn(\Delta) = -1$,~$\cyc(\Delta) = 2$ and $e(\Delta) =
xxyyxxyyyyxx$.

\begin{proposition} Let $w$ be an even word of length $2m$ and
let $s + t = 2m$. The coefficient
of $w$ in $(xy)^{\sshuf s} \shuf (yx)^{\sshuf t}$ is
\[ 2^{2m} s!t!  \sum_{\Delta}
2^{-c(\Delta)} \]
where the sum is over all the $XY$-matchings $\Delta$ such that $e(\Delta) = w$
and \hbox{$\neg(\Delta) = t$}.
\end{proposition}

\begin{proof}
By Proposition~2.2, it is sufficient to prove that
\[ N_t(w) = 2^{2m} \sum_{\Delta} %
2^{-\cyc(\Delta)} \]
where the sum is over all the $XY$-matchings $\Delta$ such that $e(\Delta) = w$
and \hbox{$\neg(\Delta) = t$}. To do this, we shall  associate to 
each $XY$-matching $\Delta = (W,\sigma)$ 
exactly $2^{2m-\cyc(\Delta)}$ $xy$-matchings
with underlying word $e(\Delta)$. 

The canonical such $xy$-matching is $\delta = (e(W),\tau)$
where $\tau$ is the involution defined by 
\[ \tau(2i-1) = 2\sigma(i) \qquad \text{if  $1 \le i \le m$} \]
For example, if $W = XYXYYX$
and $\sigma = (1234)(56)$ then 
\[ \tau =
(14)(27)(36)(58)(9\:12)(10\:11), \] as shown below.

\vbox{
\bigskip
\medskip
\bigskip

\hfil
\xymatrix@C=6pt@R=42pt{ 
X \ar@/^12pt/[rr]\ar@{-->}[d]\ar@{-->}[dr]& &
Y \ar@/_12pt/[rr]\ar@{-->}[d]\ar@{-->}[dr]& &
X \ar@/^12pt/[rr]\ar@{-->}[d]\ar@{-->}[dr]& &
Y \ar@/_15pt/@<-9pt>[llllll]\ar@{-->}[d]\ar@{-->}[dr]& & 
Y \ar@/_9pt/[rr]\ar@{-->}[d]\ar@{-->}[dr]& &
X \ar@/^15pt/[ll]\ar@{-->}[d]\ar@{-->}[dr]& \\
x\ar@{-}@/^12pt/[rrr] & 
x & 
y \ar@{-}@/_12pt/[rrr] & 
y  & 
x \ar@{-}@/^12pt/[rrr]& 
x  & 
y \ar@{-}@/_18pt/[lllll]& 
y & 
y \ar@{-}@/_15pt/[rrr] & 
y \ar@{-}@/_9pt/[r] & 
x & 
x
}
\hfill

\bigskip
\begin{center}
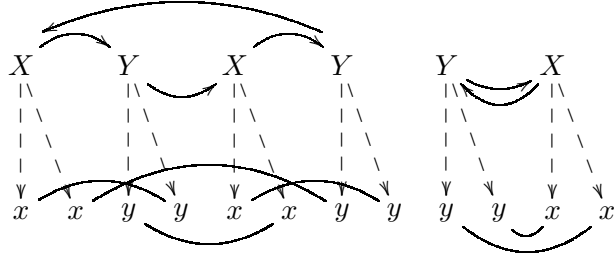
{\sc Figure 6.} The canonical $xy$-matching associated to $(W,\sigma)$.
\end{center}}
\smallskip

We obtain the remaining $xy$-matchings  by conjugating $\tau$ by
the \hbox{$2^{2m}$} 
elements of the group $\left<(12), \ldots, (2m-1\:2m)
\right>$. There is, however, some double counting, which accounts for the
factor of $2^{-\cyc(\Delta)}$. There are two cases we must consider.

Firstly, suppose $\sigma$ has a $2$-cycle, say $(k\:l)$. Then
$\tau$ involves
\[ (2k-1 \:\: 2l)(2k \:\: 2l-1), \] 
which is stabilised
by conjugation by $(2k-1 \:\: 2k)(2l-1 \:\: 2l)$. This gives us a
factor of $1/2$ for each $2$-cycle. 

Secondly, if $\hat{\sigma}$
is a cycle of length $4$ or more in $\sigma$, then we
can replace $\hat{\sigma}$ with $\hat{\sigma}^{-1}$
without changing the $xy$-matchings obtained. So again we
must compensate by a factor of $1/2$ to avoid overcounting.
For example, the two $XY$-matchings $(XYXY, (1234))$ and
$(XYXY, (1432))$ both give the same set
of $16$ $xy$-matchings.
%

Combining these observations gives the required result.
\end{proof}

\bigskip

By the same argument used to deduce Corollary~\ref{cor:nasty} 
from Proposition~\ref{prop:nasty}
we obtain the following corollary.

\begin{corollary} Let $w$ be an even word of length $2m$.
The coefficient of $w$ in $(xy-yx)^{\sshuf 2m}$ is
\[ (2m)! \vthinspace 2^{2m} \sum_{\Delta} 2^{-\cyc(\Delta)} \sgn(\Delta) \]
where the sum is over all $XY$-matchings $\Delta$ such that $e(\Delta) = w$.
Hence
\begin{equation}\label{eq:XY} 
u_{2m} = (2m)! \vthinspace 2^{2m} \sum_{\Delta} 2^{-\cyc(\Delta)} \sgn(\Delta) 
\end{equation}
where the sum is over all $XY$-matchings $\Delta$ of length $m$.
\hfill$\qed$
\end{corollary}

\subsection{} Let
\[ c_{2r} = \sum_{\Delta} \sgn(\Delta) \]
where the sum is over all $XY$-matchings $\Delta$ of length $2r$ with
\emph{just one} cycle. 
We  use an argument from the theory of exponential
structures (see~\cite[\S5.4]{StanleyII} for the general setting) to obtain
an expression for $u_{2m}$ in terms of the~$c_{2r}$.
This reduces our problem to finding the $c_{2r}$.

\begin{lemma}\label{lemma:easy}
\[ u_{2m} = (2m)!2^{2m} \sum_{\genfrac{}{}{0pt}{}{a_1,\ldots,a_m \ge 0}
{a_1 + 2a_2 + \ldots +ma_m  =m}}
\frac{(2m)!}{(2!)^{a_1} \ldots (2m)!^{a_m}} \frac{c_2^{a_1}}{2^{a_1}a_1!}
\ldots \frac{c_{2m}^{a_m}}{2^{a_m}a_m!} \]
\end{lemma}

\begin{proof}
Consider the contribution to the sum in \eqref{eq:XY} coming from 
those $XY$-matchings
whose underlying permutation has cycle type $(1^{a_1},\ldots,m^{a_m})$. 
To construct such an $XY$-matching, we must first partition $\{1 \ldots 2m\}$
into~$a_1$ subsets of size~$2$,~$a_2$ subsets of size~$4$, and so on,
up to~$a_m$ subsets of size~$2m$. This can be done
in 
\[ \frac{(2m)!}{(2!)^{a_1} \ldots (2m)!^{a_m}} \]
ways. Then we must choose for each subset an $XY$-matching with
just one cycle
on that subset. For the $a_r$ subsets chosen of size $2r$, this
can be done in~$c_{2r}^{a_r}$ ways. However we only care
about the matchings chosen, not the order we choose them in, so we must divide by~$a_r!$, giving
\[ \frac{c_2^{a_1}}{2^{a_1}a_1!}
\ldots \frac{c_{2m}^{a_m}}{2^{a_m}a_m!} \]
choices for the $XY$-matchings on the subsets. The sign of the resulting 
$XY$-matching is the product of the signs of the $XY$-matchings on the
subset, and the total number of cycles is $a_1 + \ldots + a_m$. The result
now follows from~$\eqref{eq:XY}$.
\end{proof}


\begin{lemma}
Let 
\[ f(z) = \sum_{m=0}^\infty \frac{u_{2m}}{2^{2m}(2m)!} \frac{z^{2m}}{(2m)!}
\]  be the exponential
generating function of $\frac{u_{2m}}{2^{2m}(2m)!}$. Then
\begin{equation}\label{eq:fgenfunc} 
f(z) = \exp \left( \sum_{r=1}^\infty \frac{c_{2r}}{2} \frac{z^{2r}}{(2r)!}\right). 
\end{equation}
\end{lemma}

\begin{proof}
This is an immediate consequence of Lemma~\ref{lemma:easy}. 
\end{proof}

\subsection{}
To finish the proof we need to know the $c_{2r}$.
Since the sign of an $XY$-matching is not affected by swapping the letters $X$ and $Y$, we have
\[ c_{2r} = 2\sum_{\Delta} \sgn(\Delta) \]
where the sum is over all $XY$-matchings $\Delta$ of length $2r$ such that 
$\cyc(\Delta) = 1$ and the first letter of $W$ is $X$. 

To evaluate this sum we need the following lemma.
Recall that a permutation 
$\tau$ is said to have a \emph{descent} at $i$
if $\tau(i) > \tau(i+1)$.

\begin{lemma}
Let $\Delta = (W,\sigma)$ be an $XY$-matching of length $2r$ such that $\cyc(\Delta) = 1$ and $W_1 = X$.
Suppose that $\sigma = (1 \: b_1 \ldots b_{2r-1})$.
Let $\tau \in S_{\{2 \ldots 2m\}}$ be the permutation 
defined by $\tau(i) = b_i$. Then
\[ \sgn(\Delta) = (-1)^{r-1} (-1)^{\desc(\tau)} \]
where $\desc(\tau)$ is the number of descents in $\tau$.
\end{lemma}

\begin{proof}
If  $i$ is even then $W_{b_i} = A$, and if $i$ is odd then $W_{b_i} = B$. Hence
\[ \neg(\Delta) = \#\{\text{$i$ : $i$ is even and $b_{i+1} < b_i$}\}
+  \#\{\text{$i$ : $i$ is odd and $b_{i+1} > b_i$}\}. \]
But 
$\#\{\text{$i$ : $i$ is odd and $b_{i+1} > b_i$}\} +
\#\{\text{$i$ : $i$ is odd and $b_{i+1} < b_i$}\}  = r-1$
as there are $r-1$ odd numbers $i$ such that $1 \le i < 2r-1$.
Hence
\[ \begin{split} \#\{\text{$i$ : $i$ is odd and $b_{i+1} > b_i$} \} \equiv
\#\{\text{$i$ : $i$ is odd and $b_{i+1} < b_i$}\} \\ + (r-1) \bmod 2. \end{split} \]
The result now follows.
\end{proof}

Let $\left<\genfrac{}{}{0pt}{}{t}{d}\right>$ denote the number of
permutations in $S_t$ with exactly $d$ descents.  (These
are known
as the \emph{Eulerian numbers}; our notation for them is taken from
\cite[\S5.1.3]{KnuthTAOCPIII}.) By the previous lemma,
\[ c_{2r} = 2 (-1)^{r-1} \sum_{d=0}^{2r-1} (-1)^d \left<\genfrac{}{}{0pt}{}{2r-1}{d}\right>\]
for each $r \in \N$.
By a well-known property of Eulerian numbers --- see for instance
\cite[Exercise 5.1.3(3)]{KnuthTAOCPIII} --- we have
\[  \sum_{d=0}^{2r-1} (-1)^d \left<\genfrac{}{}{0pt}{}{2r-1}{d}\right>
= (-1)^{r-1}T_r\]
where the $T_r$ are the \emph{tangent numbers}, defined by
\[ \sum_{r=1}^{\infty} T_r \frac{z^{2r-1}}{(2r-1)!} = \tan z.\]
We have therefore shown that
\begin{equation}\label{eq:g} c_{2r} = 2T_r \quad\text{for all $r\in \N$}. \end{equation}

\subsection{}
By \eqref{eq:fgenfunc} and \eqref{eq:g},
\[ \begin{split}
f(z) = \exp \left(\sum_{r=1}^\infty T_r \frac{z^{2r}}{(2r)!} \right) 
= \exp \left( \int \sum_{r=1}^\infty T_r \frac{z^{2r-1}}{(2r-1)!} \right) \\
= \exp \left( \int \tan z\right) 
 = \exp \left( \log \sec z \right) 
 = \sec z. \end{split}
\]
The right hand side is the exponential generating function for the 
Euler numbers, so
comparing coefficients, we see that $\frac{u_{2m}}{2^{2m}(2m)!} = E_{2m}$.
This completes the proof. 

\section{Appendix: Proof of Proposition~\ref{prop:nasty}}
We repeat the statement of this proposition below.

\setcounter{section}{5}
\setcounter{theorem}{2}
\begin{proposition}
Let $w$ be an even word of length $2m$ and
let $s + t = 2m$. The coefficient
of $w$ in $(xy)^{\sshuf s} \shuf (yx)^{\sshuf t}$ is
$s!t! \vthinspace N_t(w)$.
\end{proposition}

We urge the reader to read the discussion
following the original statement of this proposition in \S 5.2
before proceeding. What
follows is a formalised version of the argument we indicated there.

\begin{proof} By
an \emph{expansion} of $(xy)^{\sshuf s} \shuf (yx)^{\sshuf t}$
we mean a word of length $2m$ using
each of the letters 
\[ x^1,y^1, \ldots, x^s,y^s, y_1,x_1, \ldots 
y_t,x_t\] exactly once, and such that for each $i$, 
$x^i$ appears before $y^i$, and $x_i$ appears after $y_i$.
It should be clear that the coefficient
of $w$ in $(xy)^{\sshuf s} \shuf (yx)^{\sshuf t}$ is equal
to the number of expansions of $(xy)^{\sshuf s} \shuf (yx)^{\sshuf t}$
which become $w$ when the numbers attached to the letters are erased.

The orbits of $S_s \times S_t$ on 
expansions of $(xy)^{\sshuf s} \shuf (yx)^{\sshuf t}$
are all of size $s!t!$. Given such an orbit we shall a corresponding
$xy$-matching with underlying
word~$w$, and permutation $\sigma$. Choose any representative of
the orbit, $w^\star$ say. For $k \in \{1 \ldots 2m\}$ set $\sigma(k)= l$
where $l$ is defined by:
\[ \text{if} \:\begin{cases}
\!\text{ $w^\star_k = x_i$ then $w^\star_l = y_i$} \\
\!\text{ $w^\star_k = x^i$ then $w^\star_l = y^i$} \\
\!\text{ $w^\star_k = y_i$ then $w^\star_l = x_i$} \\
\!\text{ $w^\star_k = y^i$ then $w^\star_l = x^i$} \\
\end{cases}.
\]
It is easy to check that $\sigma$ is a fixed-point-free involution,
that $\sigma$ does not depend on the choice of $w^\star$,
and that $\neg(w,\sigma) = t$. 

Conversely, suppose we are
given an $xy$-matching with underlying word $w$ and permutation $\sigma$
and negativity $t$. Set $s = m-t$.
There is a unique
way to write $\sigma$ in the form
\[ (i^1 \: \sigma(i^1)) \ldots (i^s \: \sigma(i^s)) 
(\sigma(i_1) \: i_1) \ldots (\sigma(i_t) \: i_t) \]
such that the following conditions hold:
\begin{align*}
(1) \:\: & i^1 < i^2 < \ldots < i^s, \quad w_{i^j} = x, \quad
\sigma(i^j) > i^j 
\text{ \hskip 2pt if $1 \le j \le s$}; \\
(2) \:\: & i_1 < i_2 < \ldots < i_t\hskip 1pt, \quad w_{i_k} = x, \quad
 \sigma(i_k) < i_k 
\text{ if $1 \le k \le t$}.
\end{align*}

We shall define an associated expansion $w^\star$ of
$(xy)^{\sshuf s} \shuf (yx)^{\sshuf t}$. The underlying
word of $w^\star$ is, of course, $w$. 
For $j \in \{1 \ldots s\}$ we set $w^\star_{i^j} = x^j$ and
$w^\star_{\sigma(i^j)} = y^j$. For $k \in  \{1 \ldots t\}$
we set $w^\star_{i_j} = x_j$ and $w^\star_{\sigma(i_j)} = y_j$. 
By virtue of
our expression for $\sigma$, 
the expansion $w^\star$ we have defined is \emph{canonical}, in the sense
that the subscripts and superscripts on its letters appear in
increasing order. 

Clearly each $S_s \times S_t$-orbit contains a unique
canonical expansion.
Hence to prove the proposition
it is sufficient to prove that the two maps we have defined give a bijection
between $xy$-matchings and canonical expansions. This is merely
a matter of definition chasing. We illustrate it by an example.

Suppose $\delta = (xxyyxxyy, (13)(28)(46)(57))$. 
As $\neg(\delta) = 1$, we associate to $\delta$
an expansion of $xy \shuf xy \shuf xy \shuf yx$. Following
the given algorithm, we take $i^1 = 1, i^2 = 2, i^3 = 5$ and $i_1 = 6$
and assign the labels $(x^1 x^2 y^1 y_1 x^3 x_1 y^3 x^2)$.
Conversely, given this expansion, we get back the permutation
$(13)(28)(46)(57)$.
\end{proof}

\section*{Acknowledgements}

The first-named author wishes to thank Thierry L\'evy and Terry Lyons for
valuable discussions.


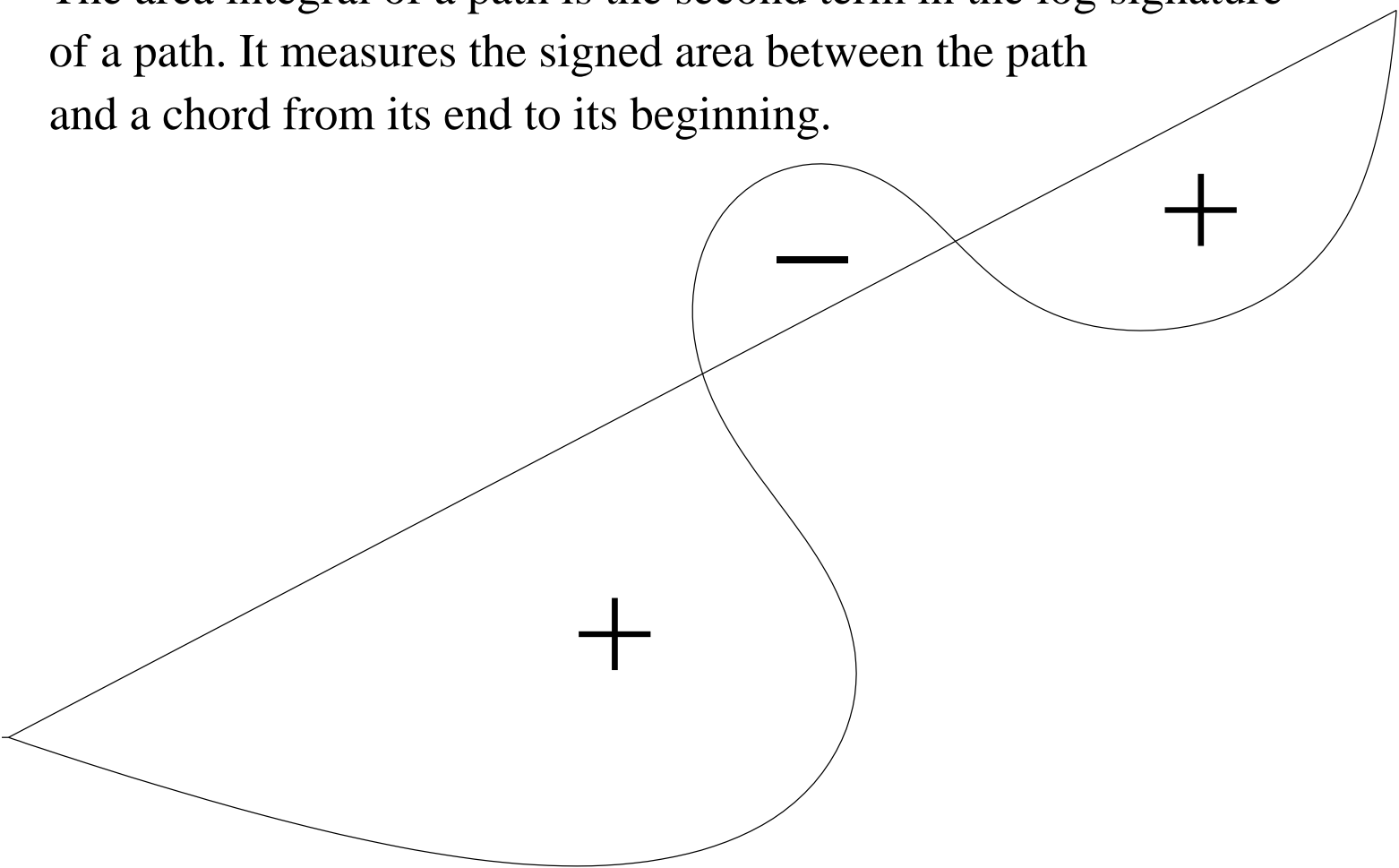
\end{document}